\documentclass{mbe}
\usepackage{amsmath}
  \usepackage{paralist}
\usepackage{bbm}
\usepackage{graphicx} 
 \usepackage[colorlinks=true]{hyperref}
\hypersetup{urlcolor=blue, citecolor=red}

\def\currentvolume{13}
 \def\currentissue{3}
  \def\currentyear{2016}
   \def\currentmonth{June}
    \def\ppages{613--619}
     \def\DOI{10.3934/mbe.2016011}

\theoremstyle{definition}

\newcommand{\pp}{\mathbb{P}}
\newcommand{\wt}{T_{\tilde{b}}}

\title[Hitting time to exponentially decaying thresholds] 
      {Approximation of the first passage time density of a Wiener process to an exponentially decaying boundary by two-piecewise linear threshold. Application to neuronal spiking activity}
\author[Massimiliano Tamborrino]{}


 \email{massimiliano.tamborrino@jku.at}

\thanks{The author is grateful to Ryota Kobayashi for the interesting discussions on adaptive-threshold models and further more.}

\begin{document}
\maketitle

\centerline{\scshape Massimiliano Tamborrino }
\medskip
{\footnotesize
 \centerline{Johannes Kepler University}
   \centerline{Altenbergerstra{\ss}e 69, 4040 Linz, Austria}
}


\begin{abstract} The first passage time density of a diffusion process to a time varying threshold is of primary interest in different fields. Here, we consider a Brownian motion in presence of an exponentially decaying threshold to model the neuronal spiking activity. Since analytical expressions of the first passage time density are not available, we propose to approximate the curved boundary by means of a continuous two-piecewise linear threshold. Explicit expressions for the first passage time density towards the new boundary are provided. First, we introduce different approximating linear thresholds. Then, we describe how to choose the optimal one minimizing the distance to the curved boundary, and hence the error in the corresponding passage time density. Theoretical means, variances and coefficients of variation given by our method are compared with empirical quantities from simulated data. Moreover, a further comparison with firing statistics derived under the assumption of a small amplitude of the time-dependent change in the threshold, is also carried out. Finally, maximum likelihood and moment estimators of the parameters of the model are derived and applied on simulated data.
\end{abstract}

\section{Introduction}
Stochastic models have been extensively used in theoretical neuroscience since the pioneer work by Gerstein and Mandelbrot in 1964 \cite{GernstenMandelbrot}. There they considered a Wiener process (also known as Brownian motion or Perfect-Integrate-and-Fire model) to model the voltage across the membrane. An action potential, also known as spike, is generated whenever the membrane potential reaches a certain constant threshold. After that, the membrane voltage is reset to its resting value and the evolution restarts. From a mathematical point of view, a spike is the first passage time (FPT) of a stochastic process to a constant threshold. The collection of spike epochs of a neuron, called spike train, defines a renewal process, with independent and identically distributed inter-spike intervals (ISIs). Despite the excellent fit with some experimental data, the Gerstein-Mandelbrot model was criticized because it disregards features involved in neuronal coding.

A first extension, combining both mathematical tractability and biological realism, is represented by \emph{Leaky-Integrate-and-Fire} (LIF) models \cite{ReviewSac,Tuckwell88}. Despite some criticisms on the lack of fit of experimental data \cite{Jolivetetal,Shinomoto}, these models are still largely used.

Another common generalization is represented by Wiener processes (or more generally LIF models) with \emph{time-dependent threshold} \cite{Tuckwell78,TuckwellWam}. These models can be chosen to reproduce biological features such as the afterhyperpolarization in neurons. For exponentially decaying thresholds, these processes can be used to model a neuron with an exponential time-dependent drift, as shown by Lindner and Longtin \cite{LindnerLongtin}. They investigated the effect of an exponentially decaying threshold on the firing statistics of a stochastic integrate-and-fire neuron \cite{LindnerLongtin}. Using a perturbation method
\cite{Lindner2004b}, they derived analytical expressions of the firing statistics under the assumption that the amplitude $\epsilon$ of the time-dependent change in the threshold is small. These statistics are useful to characterize the spontaneous neural activity and to investigate the neuronal signal transmission. In particular, they can suggest under which conditions a decaying threshold may facilitate or deteriorate signal processing by stochastic neurons. For a Wiener process, these quantities can also be obtained using the approach in \cite{Urdapilleta}. Also this method assumes a small amplitude $\epsilon$, but it has the advantage of providing an explicit approximation of the FPT density.

Here we consider a Wiener process with exponentially decaying threshold. The first aim of the paper is to provide an alternative method  to approximate the firing statistics and the FPT density for any possible amplitude $\epsilon$, extending the results in \cite{LindnerLongtin, Urdapilleta}. Different estimators are proposed, as mentioned in Section \ref{Sec1a} and discussed in Section \ref{Sec4b}. Means, variances, coefficients of variation (CVs) and distributions of the FPTs are compared on simulated data and the most suitable are recommended. A comparison with the results in \cite{LindnerLongtin, Urdapilleta} under the assumption of a small amplitude $\epsilon$ is also performed. The second aim of this work is the estimation of drift and diffusion coefficients of the Wiener process. Maximum likelihood and moment estimators are derived and evaluated on simulated data. Our results show a good approximation of both firing statistics and parameters of the underlying model.

Although the considered model generates a renewal process, the proposed method can also be applied to non-renewal processes, e.g. adaptive threshold models \cite{Chacron3,Kobayashi1}. Recently, an increasing interest arose towards these models, interest motivated by the excellent fit of the firing statistics of electrosensory neurons \cite{Chacron2,Chacron1}. The novelty of these models is that the threshold  has a jump immediately after a spike. Since the boundary depends on the previous firing epochs, the ISIs are not independent anymore. However,  the distribution between two consecutive spikes, conditioned on the initial position of the threshold, is the same of that studied here. Hence, our results may represent a first step towards an understanding of the more complicated adapting-threshold models.

\subsection{Mathematical background} \label{Sec1a}
FPTs of diffusion processes to constant or time-dependent thresholds have been extensively studied in the literature. Explicit expressions for constant thresholds are available for the Wiener process \cite{InverseGaussianBook,coxMiller}, for a special case of the Ornstein Uhlenbeck (OU) process \cite{Ricciardi}, for the Cox-Ingersoll-Ross process  \cite{CapRic}, and for those processes which can be obtained from the previous through suitable measure or space-time transformations, see e.g. \cite{Alili,CapRic,RicciardiW}.  For most of the processes arising from applications and for time-varying thresholds, analytical expressions are not available.
Numerical algorithms based on solving integral equations have been proposed in \cite{BCCP,BNR,RicciardiNip,STZ3,Taillefumier,Telve}, 
while approximations based on Monte-Carlo path-simulation methods in \cite{GS, GSZ,Metzler}.

A different approach to tackle the FPT problem consists in focusing directly on the two-sided boundary crossing probability (BCP), i.e. the probability that a process is constrained to be between two boundaries. If one of the boundary is set to $-\infty$, the resulting one-sided BCP equals the survival probability of the FPT to the other boundary \cite{Wang1997}. Explicit formulas for the BCP of a standard Wiener process for continuous and piecewise-linear thresholds are known (see \cite{BorovkovNovikov,Novikovetal,Wang2001,Wang1997,Wang2007}). In general, the BCP of a diffusion process through an exponential decaying threshold is available only for those processes which can be expressed as a piecewise monotone functional of a standard Brownian motion. Examples are the OU process or the geometric Brownian motion with time dependent drift for specific parameter values \cite{Wang2007}.
The simple but powerful idea is to approximate both one and two-sided curvilinear BCPs by similar probabilities for close boundaries of simpler form, namely $n$ piecewise-linear thresholds, whose computation of the BCP for Wiener is feasible. Under some mild assumptions, the approximated two-sided BCP converges to the original one when $n\to\infty$ \cite{Wang2007}, with rate of convergence given in \cite{BorovkovNovikov,Wang1997}.

For the exponential decaying threshold considered in this paper,
the convergence can be obtained by choosing piecewise linear thresholds approximating the curved boundary from above and below, with approximation accuracy given by their distance \cite{Wang2007}. However, all the available formulas for the BCPs require either Monte-Carlo simulation methods or heavy numerical approximations.

Here we consider a two-piecewise linear threshold as an approximation of the  curvilinear boundary. Since $n=2$, the asymptotic convergence of the BCPs does not hold. However, we can derive analytical expression for the FPT density to the two-piecewise linear boundary, and use it as an approximation of the unknown FPT density. Four possible piecewise thresholds are proposed and optimized to minimize the distance to the original threshold.

\section{Model}\label{Sec2}
\begin{figure}
\includegraphics[width=\textwidth]{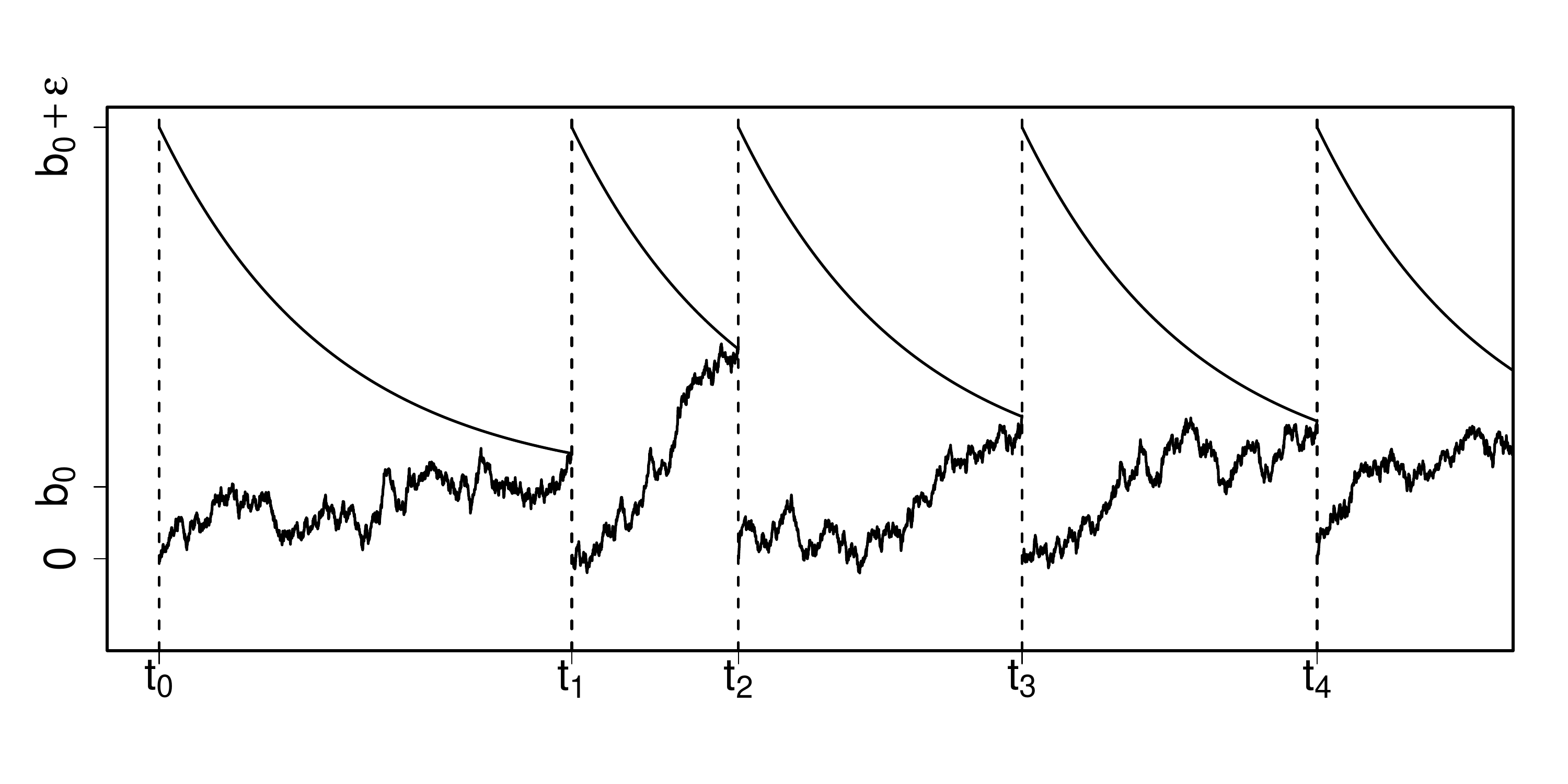}
\caption{Schematic illustration of the single trial of a Wiener process in presence of exponentially decaying threshold $b(t)=b_0+\epsilon\exp(-\lambda (t-\delta_k))$, where $\delta_k$ denotes the $k$th spike. The membrane potential starts in $x_0=0$ at time $\delta_0:=t_0=0$, and it evolves until it hits the boundary for the first time. Then, a spike is generated, the voltage $X(t)$ is reset to its resting potential $x_0$, the threshold is reset to $b_0+\epsilon$ and the evolution restarts. For the considered spiking generation rule, all ISIs are independent and identically distributed. Here the parameters are $\mu=1, \sigma^2=1, b_0=1, \lambda=1$ and $\epsilon=5$.}
\label{FigFPT}
\end{figure}

We describe the membrane potential evolution of a single neuron by a Wiener process $X(t)$, starting at some initial value $x_0$. We assume $X(t)$ given as the solution to a stochastic differential equation
\begin{equation}\label{model}
\left\{
\begin{array}{l}
dX(t)=\mu dt + \sigma dW(t),\\
X(t_0)=x_0,  \qquad t>t_0,
\end{array}
\right.
\end{equation}
where $W(t)$ is a standard (driftless) Wiener process. The drift $\mu>0$ and the diffusion coefficient $\sigma>0$ represent input and noise intensities, respectively. A spike occurs when the membrane potential $X(t)$ exceeds the exponentially decaying threshold
\begin{equation}\label{b}
b^*(t)=b_0+\epsilon \exp\left[-\lambda (t-\delta_k)\right]
\end{equation}
for the first time. Here, $\delta_k$ denotes the time of the $k$th spike for $k>0$, and can be interpreted as a relative refractory period. We set $\delta_0$ to be the starting time of the process, i.e. $\delta_0=t_0$. The term $\lambda$ represents the decay rate of the threshold, while $\epsilon$ is interpreted as the amplitude of the time-dependent change in the boundary. After a spike, the membrane potential is reset to its resting position $x_0<b_0+\epsilon$, and its evolution is restarted, as illustrated in Fig. \ref{FigFPT}. The presence of $\delta_k$ in \eqref{b} ensures that the ISIs are independent and identically distributed. Denote by $b(t)$ the threshold $b^*(t)$ for $k=0$, i.e.
\begin{equation*}
b(t)=b_0+\epsilon \exp\left[-\lambda (t-t_0)\right].
\end{equation*}
Then, all ISIs are distributed as the FPT of $X$ to $b(t)$, namely
\begin{equation*}
T_b=\inf\{t>t_0: X(t)\geq b(t)\}.
\end{equation*}
Quantities of interest are the probability density function (pdf) and the cumulative distribution function (cdf) of $T_b$, denoted by $f_{T_b}$ and $F_{T_b}$, respectively. Another relevant quantity is the two-sided BCP given by
\begin{equation*}
\pp_X(a,c,\tau)=\pp\left(a(t)<X(t)<c(t), \forall t \in [t_0,\tau]\right).
\end{equation*}
Here $\tau>t_0$ is fixed, boundaries $a(t)$ and $c(t)$ are real functions satisfying $a(t)<c(t)$ for all $t_0<t\leq \tau$ and $a(t_0)<x_0<c(t_0)$.  Setting $a(t)=-\infty$ yields the one-sided BCP
\[
\pp_X(-\infty,c,\tau)=\pp(T_c>\tau)=1-F_{T_c}(\tau),
\]
which corresponds to the survival probability of $T_c$. For a standard Wiener process $W$, Wang and P\"{o}tzelberger \cite{Wang2007} showed that, if the sequences of piecewise linear functions $a_n$ and $c_n$ converge uniformly to $a(t)$ and $c(t)$ on $[t_0,\tau]$ respectively, then, for the continuity property of probability measure, it holds
\begin{equation*}
\lim_{n\to \infty} P_W(a_n,c_n,\tau)=P_W(a,c,\tau).
\end{equation*}
When $a(t)=-\infty$ and $c(t)=b(t)$, the convergence of $\pp(T_{b_n}>\tau)$ to $\pp(T_b>\tau)$ can be obtained by choosing piecewise linear thresholds approximating $b(t)$ from above, denoted by  $b_n^+(t)$, or from below, $b_n^-(t)$. That is, $b_n^+(t)\geq b(t)$ and $b_n^-(t)\leq b(t), \ \forall t\in[t_0,\tau]$, respectively. Since the considered curved boundary is convex, we have \cite{Wang2001}
\begin{equation}\label{updown}
\pp_X(-\infty,b_n^-,\tau)\leq \pp_X(-\infty,b,\tau)\leq \pp(-\infty,b_n^+,\tau),
\end{equation}
i.e.
\begin{equation*}
\pp(T_{b^+_n}\leq \tau)\leq \pp(T_b\leq \tau)\leq \pp(T_{b^-_n}\leq \tau).
\end{equation*}
The approximation accuracy is given by $\pp_X(-\infty,b_n^+,\tau)-\pp_X(-\infty,b_n^-,\tau)=F_{T_{b^-_n}}(\tau)-F_{T_{b^+_n}}(\tau)$, with bounds given in \cite{BorovkovNovikov}. Obviously, the accuracy in the BCP increases when the distance between the two thresholds decreases.

\section{FPT to continuous piecewise linear threshold}\label{Method}
The transition density function of a standard Brownian motion in $x_1, x_2,\ldots , x_n$ at time $t_1, t_2, \ldots, t_n$, constrained to be below the absorbing threshold $c(t)$ defined by $n$ piecewise-linear threshold over $[t_0,t_n]$, is given in \cite{Wang1997}. Extending that result to the case of a Brownian motion with drift $\mu$ and diffusion coefficient $\sigma$, starting in $x_0<c(t_0)=c_0$ at time $t_0$, we obtain
\begin{eqnarray}\label{eq2.9}
\nonumber&&p_{c}(x_1,t_1;x_2,t_2;\ldots;x_n,t_n)=\prod_{i=1}^n p_{c}(x_i,t_i|x_{i-1},t_{i-1})\\
&=&\prod_{i=1}^n \frac{\left[1-\exp\left(-\frac{2(c_i-x_i)(c_{i-1}-x_{i-1})}{\sigma^2(t_i-t_{i-1})}\right)\right]}{\sqrt{2\pi\sigma^2(t_i-t_{i-1})}}\exp\left(-\frac{[x_i-x_{i-1}-\mu(t_i-t_{i-1})]^2}{2\sigma^2(t_i-t_{i-1})}\right),\quad  
\end{eqnarray}
where $c_i=c(t_i)$ and $x_i< c_i, 1\leq i\leq n$.

\noindent From \eqref{eq2.9}, it follows that \cite{ZuccaSac}
{\small{\begin{align}\label{eq2.10}
\pp(W_{t_1}\in C_1,\ldots,W_{t_n}\in C_n,T_c>t_n)=\int_{C_1}\cdots\int_{C_n}p_{c}(x_1,t_1;\ldots;x_n,t_n|x_0,t_0)dx_1\cdots dx_n,
\end{align}}}
for any Borel set $C_i\subseteq (-\infty,c_i), 1\leq i \leq n$. If $C_i=(-\infty,c_i)$, then \eqref{eq2.10} is equal to $\pp(T_c>t_n)$, and it holds
\begin{equation}\label{fpt}
f_{T_c}(t)=-\frac{\partial }{\partial t_n} \int_{-\infty}^{c_1}\cdots\int_{-\infty}^{c_n} p_{c}(x_1,t_1;\cdots,x_n,t_n|x_0,t_0) dx_1\cdots dx_n.
\end{equation}

When $n=1$, the pdf $f_{T_c}$ is known. Since $X$ is a Wiener process with positive drift, the distribution of the FPT to $c(t)=\alpha+\beta (t-t_0)$ is inverse Gaussian,
$T_{c}\sim IG\left[(\alpha-x_0)/(\mu-\beta), (\alpha-x_0)^2/\sigma^2\right]$, with pdf
\begin{equation}\label{fpt3}
f_{T_c}(t)=\frac{\alpha-x_0}{\sqrt{2\pi\sigma^2(t-t_0)^3}}\exp\left(-\frac{\left[\alpha-x_0-(\mu-\beta)(t-t_0)\right]^2}{2\sigma^2(t-t_0)}\right),
\end{equation}
mean $\mathbb{E}[T_{c}]=(\alpha-x_0)/(\mu-\beta)$ and variance $\textrm{Var}(T_{c})=(\alpha-x_0)\sigma^2/(\mu-\beta)^3$ \cite{InverseGaussianBook,coxMiller}. Note that the distribution of $T_c$ is the same of that of the FPT of a Wiener process with positive drift $\mu-\beta$ to a constant threshold $c(t)=\alpha$.  In general, the approximation of $F_{T_b}$ by $F_{T_c}$ when $n=1$ is too rough. However, when $\lambda$ is very small, $\exp(-\lambda t)\approx 1-\lambda t$, yielding $b(t)\approx b_0+\epsilon-\lambda \epsilon t$. Hence, $T_b$ can be approximated by $T_c$ with $\alpha=b_0+\epsilon$ and $\beta=-\lambda \epsilon$.

Since we approximate $b(t)$ by means of a continuous two-piecewise linear threshold, we denote by $\tilde b$ the linear threshold $c(t)$ when $n=2$. We have
\begin{equation}\label{St}
\tilde{b}(t)=\tilde b_1(t)\mathbbm{1}_{\{t\leq t_1\}}+\tilde b_2(t) \mathbbm{1}_{\{t>t_1\}}=\left\{
\begin{array}{ll}
\alpha_1+\beta_1 (t-t_0) & \textrm{ if } t_0\leq t\leq t_1\\
\alpha_2+\beta_2(t-t_1) & \textrm{ if } t> t_1
\end{array}
\right. ,
\end{equation}
where $\mathbbm{1}_A$ denotes the indicator function of the set $A$ and $\alpha_1,\alpha_2,\beta_1,\beta_2\in \mathbb{R}$. Throughout the paper, we set $\alpha_2=\alpha_1+\beta_1 (t_1-t_0)$ to guarantee the continuity of $\tilde b(t)$. This allows to provide analytical expressions of \eqref{eq2.9} and \eqref{fpt}, which we use as an approximation of $f_{T_b}$. In particular, we have
\begin{align}\label{fpt2}
\nonumber\pp(\wt<t)=&\ \pp(\wt< t,\wt<t_1)+\pp(\wt<t,\wt>t_1)\\
\nonumber=&\  \pp(T_{\tilde{b}_1} <\min(t_1,t)) + \int_{-\infty}^{\tilde{b}(t_1)} \pp(T_{\tilde{b}_2}<t|X(t_1)=x_1) p_{\tilde{b}_1}(x_1,t_1) dx_1\\
=&  \int_0^{\min(t,t_1)} f_{T_{\tilde{b}_1}}(s) ds +
\int_{-\infty}^{\tilde{b}(t_1)} \int_{t_1}^t f_{T_{\tilde{b}_2}}(s|x_1,t_1) p_{\tilde{b}_1}(x_1,t_1)  ds dx_1 
\end{align}
with $p_{\tilde{b}_1}(x_1,t_1)$ given by \eqref{eq2.9} for $n=1$. Mimicking \cite{TDL3}, by taking the derivative of \eqref{fpt2} with respect to $t$, and plugging \eqref{fpt3} in it for proper values of $\alpha$ and $\beta$, we get
{\small\begin{align}\label{fpt4} 
\nonumber &\hspace{-.3cm} f_{\wt}(t)\\ \nonumber
=&f_{IG\left(\frac{\alpha_1-x_0}{\mu-\beta_1},\frac{(\alpha_1-x_0)^2}{\sigma^2}\right)}(t-t_0)\mathbbm{1}_{\{t\leq t_1\}}+\int_{-\infty}^{\tilde{b}(t_1)} f_{IG\left(\frac{\alpha_2-x_1}{\mu-\beta_2},\frac{(\alpha_2-x_1)^2}{\sigma^2}\right)}(t-t_1)p_{\tilde{b}_1}(x_1,t_1) dx_1\\
\nonumber =& \frac{\alpha_1-x_0}{\sqrt{2\pi\sigma^2(t-t_0)^3}} \exp\left(-\frac{(\alpha_1-x_0-(\mu-\beta_1)(t-t_0))^2}{2\sigma^2(t-t_0)}\right)\mathbbm{1}_{\{t\leq t_1\}}\\
\nonumber&+\mathbbm{1}_{\{t> t_1\}}
\frac{1}{\sqrt{2\pi\sigma^2(t-t_0)^3}}\exp\left(-\frac{(\alpha_2-x_0-(\mu-\beta_2)(t-t_1)-\mu(t_1-t_0))^2}{2\sigma^2(t-t_0)}\right)\\
\nonumber&\times \left\{[\alpha_2-x_0-\beta_2(t_1-t_0)]\Phi\left(\frac{(\alpha_2-x_0-\beta_2(t_1-t_0))\sqrt{(t-t_1)}}{\sqrt{\sigma^2(t_1-t_0)(t-t_0)}}\right)\right.\\\nonumber& \left.-(\alpha_2+x_0-\beta_2(t_1-t_0)-2\alpha_1)
\exp\left(-\frac{2(t-t_1)(\alpha_1-x_0)(\alpha_2-\alpha_1-\beta_2(t_1-t_0))}{\sigma^2(t_1-t_0)(t-t_0)}\right)\right.\\& \times\left.\Phi\left(\frac{(\alpha_2+x_0-\beta_2(t_1-t_0)-2\alpha_1)\sqrt{(t-t_1)}}{\sqrt{\sigma^2(t_1-t_0)(t-t_0)}}\right)\right\}.
\end{align}}This result extends that for a driftless Brownian motion, see e.g. \cite{Abundo,Scheike}.
%
As expected, setting $\alpha_1=\alpha_2=\alpha$ and $\beta_1=\beta_2=\beta$ yields the pdf of the FPT of a Wiener process to a linear  threshold $c(t)=\alpha+\beta(t-t_0)$. By definition, the first two moments and variance of $T_{\tilde b}$ are given by
\begin{equation}\label{EVT}
\mathbb{E}[T_{\tilde b}]=\int_0^\infty tf_{T_{\tilde b}}(t)dt, \qquad
\mathbb{E}[T^2_{\tilde b}]=\int_0^\infty t^2 f_{T_{\tilde b}}(t)dt, \qquad \textrm{Var}[T_{\tilde b}]=\mathbb{E}[T_{\tilde b}^2]-\mathbb{E}[T_{\tilde b}]^2,
\end{equation}
and can be numerically computed.

\section{Parameter estimation}\label{Sec4}
\subsection{Parameter estimation of the piecewise-linear threshold}\label{Sec4a}
The primary aim of this paper is the approximation of the FPT distribution (and relevant statistics) for a curved boundary $b(t)$, by means of the FPT distribution for a continuous two-piecewise linear threshold $\tilde b(t)$. As discussed in Section \ref{Sec2}, the quality of the approximation improves when the distance between $\tilde b$ and $b$ decreases. 

\begin{figure}
\includegraphics[width=1.0\textwidth]{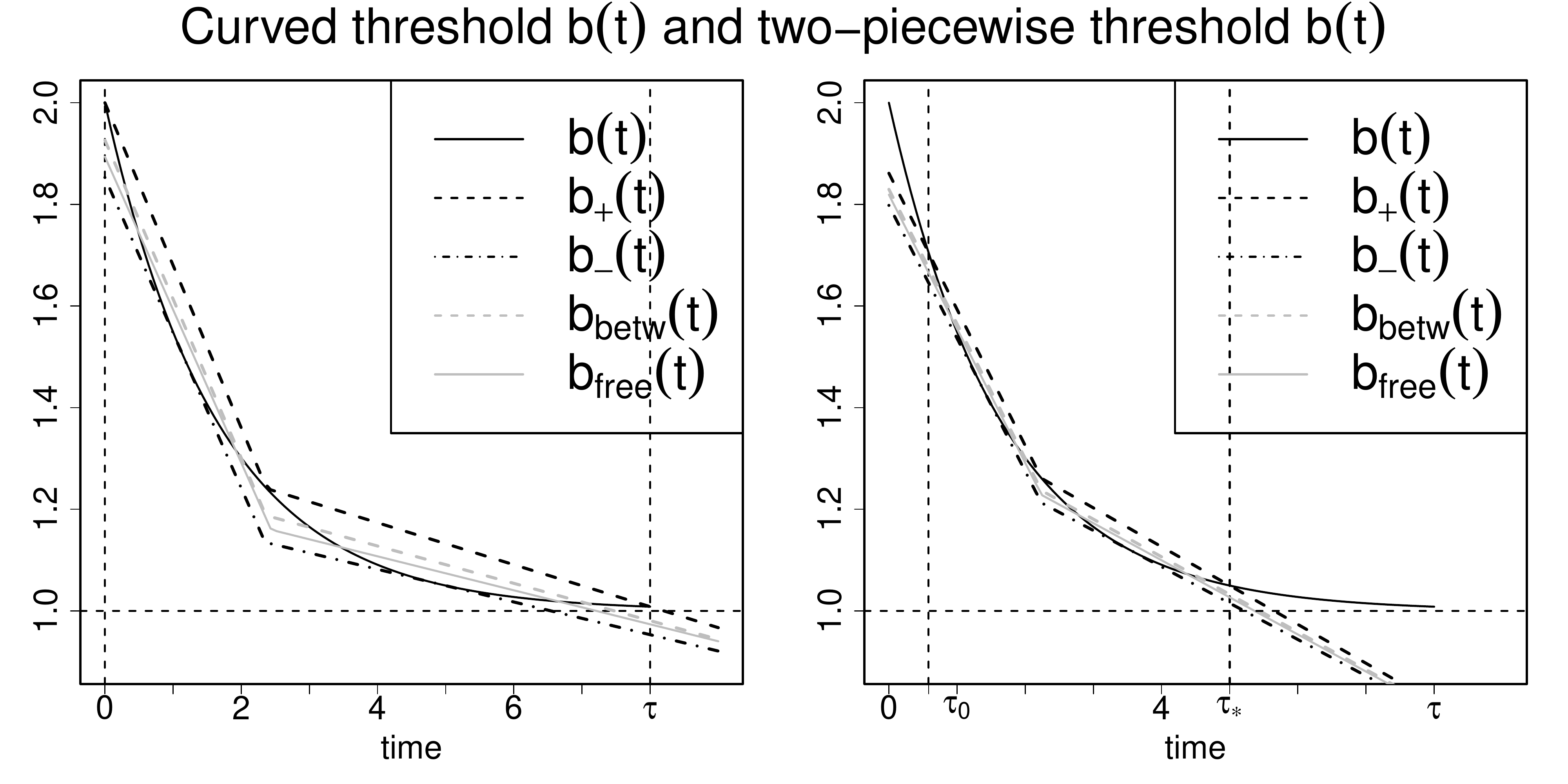}
\caption{Curved threshold $b(t)$ (continuous line) and four proposed approximating piecewise-linear thresholds: $b_+(t)$ from above (dashed lines) ; $b_-(t)$ from below (dashed-dotted line); $b_\textrm{betw}(t)$ which is equidistant from $b_+$ and $b_-$ (gray dashed line); $b_\textrm{free}(t)$ with no restrictions (gray continuous line). For each type of linear threshold, the best approximation is given by the line minimizing a function of the distance from $b$ on $[t_0=0,\tau]$ (left figure) and on $[\tau_0,\tau_*]\subseteq [t_0=0,\tau]$ (right figure). As discussed in Section \ref{Sec4a}, a shorter time interval provides a better approximation of $b$.}
\label{thresh}
\end{figure}

Denote by $\theta=(\alpha_1,\beta_1,\beta_2,t_1)$ the parameters of $\tilde b$ in \eqref{St}, with $\alpha_2=\alpha_1+\beta_1(t_1-t_0)$. We are interested in determining the estimator $\hat\theta$  which minimizes $|\tilde b(t)-b(t)|$ on $[\tau_0,\tau_*] $, with $t_0<\tau_0<t_1<\tau_*<\tau$. The time interval is chosen such that the probability of having a FPT outside it is smaller than $0.005$, i.e.
\begin{equation}\label{mah}
\pp(T_b \in [\tau_0,\tau_*])\geq 0.99.
\end{equation}
Doing this, we improve the approximation of $b$ on $[\tau_0,\tau_*]$ (cf. Fig. \ref{thresh}), allowing a larger deviation from $b$ on $[t_0,\tau_0)$ and $(\tau_*,\tau]$, i.e. on intervals where the probability of observing a FPT is low. Since
\begin{eqnarray*}
\pp(T_b>t)=\pp(X(s)<b(s), s \in [0,t])\leq \pp(X(t)<b(t))
\end{eqnarray*}
and $b(t)>b_0$, it follows that
\[
\pp(X(t)\geq b(t)) \leq \pp(T_b\leq t)\leq \pp(T_{b_0} \leq t),
\]
with $T_{b_0}\sim IG((b_0-x_0)/\mu, (b_0-x_0)^2/\sigma^2)$.  Since $X$ is a Wiener process, $X(t)\sim N(\mu t, \sigma^2 t)$. Then, we choose $\tau_0$ and $\tau_*$ such that
\[
\pp(T_{b_0}\leq\tau_0)=0.005, \qquad
\pp(X(\tau_*)\geq b(\tau_*))=0.995,
\]
yielding the desired probability \eqref{mah}.

Throughout the paper, we consider four possible continuous two-piecewise linear boundaries on $[\tau_0,\tau_*]$, as illustrated in Fig. \ref{thresh}:
\begin{enumerate}
\item Threshold $b_+$ approximating  $b$ from above on $[\tau_0,\tau_*]$, passing through $(\tau_0,\break b(\tau_0))$, $(t_1,b(t_1))$ and
 $(\tau_*,b(\tau_*))$,
\begin{equation*}\label{th2}
b_+(t)=b(\tau_0)+\frac{b(t_1)-b(\tau_0)}{t_1-\tau_0}(t-\tau_0)\mathbbm{1}_{\{t\leq t_1\}}+
\frac{b(\tau_*)-b(t_1)}{\tau_*-t_1}(t-t_1)\mathbbm{1}_{\{t>t_1\}},
 \end{equation*}
 i.e.
 \[
 \alpha_1=b_+(t_0), \quad \beta_1= \frac{b(t_1)-b(\tau_0)}{t_1-\tau_0}, \quad \beta_2=\frac{b(\tau_*)-b(t_1)}{\tau_*-t_1}(t-t_1).
 \]
Due to the assumptions, for given $\tau_0$ and $\tau_*$, $t_1$ is the only unknown quantity.
\item Threshold $b_-$ approximating $b$ from below on $[\tau_0,\tau_*]$. We assume that $b_-$ is tangent to $b(t)$ in both $\tilde t_1$ and  $\tilde t_2>\tilde t_1$, with $t_1$ intersection time point of the two tangent lines
\[
y_i(t)=b(\tilde t_i)-\lambda \epsilon \exp(-\lambda \tilde t_i)(t-\tilde t_i),
\]
for $i=1,2$. Setting $y_1(t_1)=y_2(t_1)$, we get
\[
t_1=\frac{\exp(-\lambda \tilde t_1)[1+\lambda \tilde t_1]-\exp(-\lambda \tilde t_2)[1+\lambda \tilde t_2]}{\lambda[\exp(-\lambda \tilde t_1)-\exp(-\lambda \tilde t_2)]}.
\]
Then, the desired threshold $b_-(t)$ is
\[
b_-(t)=y_1(\tilde t_1)+\frac{y_1(t_1)-y_1(\tilde t_1)}{t_1-\tilde t_1}(t-\tilde t_1)\mathbbm{1}_{\{t\leq t_1\}}+\frac{y_2(\tilde t_2)-y_2(t_1)}{\tilde t_2-t_1}(t- t_1)\mathbbm{1}_{\{t>t_1\}},
\]
with
\[
\alpha_1=b_-(t_0),\quad \beta_1=\frac{y_1(t_1)-y_1(\tilde t_1)}{t_1-\tilde t_1}, \quad
\beta_2=\frac{y_2(\tilde t_2)-y_2(t_1)}{\tilde t_2-t_1}.
\]
For fixed $\tau_0$ and $\tau_*$, the unknown parameters are $\tilde t_1$ and $\tilde t_2$.
\item Threshold $b_\textrm{betw}(t)$ constrained to be between $b_+(t)$ and $b_-(t)$ on $[\tau_0,\tau_*]$, i.e.\break $b_-(t)\leq b_\textrm{betw}(t)\leq b_+(t)$.
\item Threshold with no constraints, denoted by $b_\textrm{free}(t)$.
\end{enumerate}
Denote by $\hat\theta_+, \hat\theta_-,\hat\theta_\textrm{betw}$ and $\hat\theta_\textrm{free}$ the estimators of $\theta$ from the boundaries $b_+, b_-,b_\textrm{betw}$ and $b_\textrm{free}$, respectively.
From \eqref{updown}, it follows that the best approximation of $\pp_X(-\infty,b,\break\tau)$ is obtained when the distance between $b_+$ and $b_-$ is minimized. For this reason, we define $\hat\theta_+$ and $\hat\theta_-$ as the estimators minimizing the area of the squared distance between the two boundaries, i.e.
\begin{equation*}
(\hat\theta_+,\hat\theta_-)=\arg\min_{(\theta_+,\theta_-)}\left[ \int_{\tau_0}^{\tau{^*}} |b_+(t)-b_-(t)|^2dt\right],
\end{equation*}
with $\hat\theta_+$ and $\hat\theta_-$ satisfying the conditions $b_+(t)>b(t)$ and $b_-(t)<b(t)$ on $[\tau_0,\tau_*]$. The quantity $|b_+(t)-b_-(t)|^2$ instead of $|b_+(t)-b_-(t)|$ is chosen to avoid possible numerical issues in the optimization procedure.

Once $b_+$ and $b_-$ have been computed, the estimator $\hat\theta_\textrm{betw}$ is defined as
\[
\hat\theta_\textrm{betw}=\arg\min_{\theta}\left[ \int_{\tau_0}^{\tau{^*}} \left(|b_+(t)-b_\textrm{betw}(t)|^2 + |b_-(t)-b_\textrm{betw}(t)|^2\right)dt\right],
\]
i.e. it is the equidistant line from $b_+(t)$ and $b_-(t)$.

Finally, the estimator $\hat\theta_\textrm{free}$ is the one minimizing the area of the squared distance between the piecewise-line and the curved boundary, i.e.
\[
\hat\theta_\textrm{free}=\arg\min_{\theta}\left[\int_{\tau_0}^{\tau_*} |b_\textrm{free}(t)-b(t)|^2 dt\right].
\]
Note that the estimation of $\theta$ \emph{does not depend} on the observations of the FPTs but it is performed theoretically under the assumption that the parameters $\lambda, \epsilon$ and $b_0$ of $b(t)$ are known.
The proposed estimators and their assumptions are summarized in Table \ref{Table1}. All the minimizations have been performed in the computing environment \textbf{R}  \cite{R}. Since the parameter values need to fulfil some conditions (cf. Table \ref{Table1}), minimizing the areas is a constrained optimization problem. We perform it by means of the built-in \textbf{R} function \verb|optim|, penalizing those parameter values not fulfilling the conditions by returning $10^{10}$.

\begin{table}[b]
\scalebox{0.9}{\begin{tabular}{|c|c|c|c|}
\hline
Estimator& Assumption on $\tilde b(t)$ on $[\tau_0,\tau_*]$ & Unknown parameters & Parameter conditions \\
\hline
$\hat\theta_+$& $b_+(t)\geq b(t)$& $t_1$& $t_1\in[\tau_0,\tau_*]$ \\
\hline
$\hat\theta_-$& $b_-(t)\leq b(t)$ & $\tilde t_1,\tilde t_2$ & $\tau_0\leq \tilde t_1\leq \tilde t_2\leq \tau_*$\\
\hline
$\hat\theta_\textrm{betw}$& $b_-(t)\leq b_\textrm{betw}(t)\leq b_+(t)$ & none& none \\
\hline
$\hat\theta_\textrm{free} $& none & $\alpha_1,\beta_1,\beta_2,t_1$ & none \\
\hline
\end{tabular}}
\caption{Proposed estimators of the parameters of the piecewise linear thresholds $b_+,b_-,b_\textrm{betw}$ and $b_\textrm{free}$ given in Section \ref{Sec4a} under different assumptions.}
\label{Table1}
\end{table}

\subsection{Parameter estimation of the process} \label{Sec4b}
The second aim of the paper is the estimation of the parameters $\mu$ and $\sigma^2$ of the Wiener process  from a sample $\{r_i\}_{i=1}^n$ of $n$ independent observations of $T_b$. That is, we want to estimate $\phi=(\mu,\sigma^2)$ under the assumption that the parameters of the threshold are known.

\subsubsection{Maximum likelihood estimator of $\phi$} First, we derive the parameters $\theta$ of the threshold $\tilde b(t)$, as described in Section \ref{Sec4a}. Then, we use maximum likelihood estimator (MLE) as follows. Since the observations are independent and identically distributed, the log-likelihood function is given by
\begin{equation*}
l_r(\phi)=\sum_{i=1}^n \log f_{T_{\tilde{b}}}(r_i;\phi),
\end{equation*}
where $f_{T_{\tilde b}}$ is the pdf given by \eqref{fpt4} with $\theta$ replaced by $\hat\theta$. Then, the  log-likelihood function can be maximized numerically to obtain the unknown parameter $\phi$. Since the parameter values of $\mu$ and $\sigma$ need to be positive, when minimizing $l_r(\phi)$ by means of the function \verb|optim|, we penalize negative values of $\mu$ and $\sigma$ by returning $10^{10}$. We denote by $\hat\phi_{\textrm{MLE}}(\hat\theta)$ the MLE of $\phi$ derived from the threshold $\tilde b$ with parameters $\hat\theta$.

\subsubsection{Moment estimator} A different approach consists in equating the theoretical moments of $T_{\tilde b}$, given by Eq. \eqref{EVT}, with the empirical moments of $T_b$. In particular, we numerically solve a system of two equations (given by the first two moments) in the two unknown parameters $\phi=(\mu,\sigma^2)$. We denote by $\hat\phi_{\textrm{ME}}$ the moment estimator (ME) of $\phi$.

When $\epsilon$ is small, approximated mean and variance of $T_b$  are available \cite{LindnerLongtin,Urdapilleta}. In particular, for a general parameter $b_0>x_0$, we have
\begin{eqnarray}
\label{ET} \widehat{\mathbb{E}[T_b]}&=& \frac{b_0}{\mu}+\frac{\epsilon}{\mu}\exp\left(\frac{b_0\left(\mu-\sqrt{\mu^2+2\lambda\sigma^2}\right)}{\sigma^2}\right),\\
\nonumber \widehat{\textrm{Var}(T_b)}&=&\frac{b_0\sigma^2}{\mu^3}+\frac{\sigma^2\epsilon}{\mu^3}(b_0-1)\\
&+&\frac{2\epsilon}{\mu^2}\left(\frac{\mu b_0}{\sqrt{\mu^2+2\lambda\sigma^2}}+\frac{\sigma^2}{2\mu}-\theta_0\right)\exp\left(\frac{b_0\left(\mu-\sqrt{\mu^2+2\lambda\sigma^2}\right)}{\sigma^2}\right).\quad
\label{VT}
\end{eqnarray}
We denote by $\hat\phi_{\textrm{ME}}^\epsilon$ the moment estimator of $\phi$ obtained from  \eqref{ET} and \eqref{VT} when $\epsilon$ is small.

\section{Simulation study}

\subsection{Monte Carlo simulations} We simulate FPTs of the Wiener process $X$ to $b(t)$ as described in \cite{LindnerLongtin,ReviewSac}. Applying the Euler-Maruyama scheme to the stochastic differential equation \eqref{model}, we generate realizations of $X$, denoted by  $x_i:=X(s_i)$, at discrete times $s_i=i \Delta s, i\geq1$. We set $X_0=x_0=0$ and $\Delta s=0.001$ as time step. To avoid the risk of not detecting a crossing of the boundary due to the discretization of the sample path, at each iteration step we compute the probability that the bridge process $X^{[s_i,s_{i+1}]}=\left\{X_s^{[s_i,s_{i+1}]},s\in [s_i,s_{i+1}]\right\}$, originated in $x_i<b(s_i)$ at time $s_i$ and constrained to be in $x_{i+1}<b(s_{i+1})$ at time $s_{i+1}$, crosses the threshold in between $s_i$ and $s_{i+1}$. For a Wiener process, this probability is given by \cite{Honerkamp}
\begin{equation*}
\pp(x_i,x_{i+1})=\exp\left\{-\frac{2[b(s_{i+1})^2-b(s_{i+1})(x_i+x_{i+1})+x_ix_{i+1}]}{\sigma^2 \Delta s}\right\}.
\end{equation*}
A FPT is observed if $x_i$ hits/exceeds the threshold $b$ at time $s_i$, i.e. $x_i\geq b(s_i)$, or if the probability of having crossed the threshold in $(s_i,s_{i+1})$ is larger than a randomly generated uniform number $u_i$ in $(0,1)$, i.e. $\pp(x_i,x_{i+1})>u_i$. In this case, the mid-point $(s_i+s_{i+1})/2$ is chosen as simulated FPT. Samples of size $100$ are simulated for different values of $\sigma^2, \lambda$ and $\epsilon$ when $b_0=1$ and $\mu=1$. In particular, we consider $\sigma^2=0.2, 0.4, 1$; $\epsilon=0.05, 0.1, 0.2, 1, 5, 10$ and $\lambda=0.02, 0.04, 0.08, 0.15, 0.30, 0.60, 1.00, 3.00, 5.00, 10.00$. These parameter values are chosen to cover and extend the cases of small values of $\epsilon$ considered in \cite{LindnerLongtin,Urdapilleta}. Finally, for each value of $\sigma^2,\epsilon$ and $\lambda$, we repeat simulation of data set 1000 times, obtaining 1000 statistically indistinguishable and independent trials.

\subsection{Set up} In the simulations we are mainly concerned to illustrate the performance of our method under the assumption that the threshold $b(t)$ is known, i.e. $b_0$, the rate $\lambda$ and the amplitude $\epsilon$ are given. Two scenarios are considered: both $\mu$ and $\sigma^2$ are known; no information about the parameter of the Wiener process is given.
In the first case it is of interest to evaluate the error in the estimation of mean, variance, CV and cdf of $T_b$ by comparing  theoretical \eqref{EVT} and empirical firing statistics. When $\epsilon$ is small, a further comparison with \eqref{ET} and \eqref{VT} is  carried out. To measure the error in the estimation of $F_{T_b}$, we consider the relative integrate absolute error ($R_\textrm{IAE}$), defined as
\begin{equation}\label{riae}
R_{\textrm{IAE}}(\hat F_{T_b})=\frac{\int_0^{\infty} | \hat F_{T_b}(t)-F_{T_b}(t)| dt}{\mathbb{E}[T_b]}.
\end{equation}
We replace the unknown quantities $F_{T_b}$ and $\mathbb{E}[T_b]$ by their empirical estimators, defined by $F_n(t)=\frac{1}{n}\sum_{i=1}^n \mathbbm{1}_{\{T_{b_i}\leq t\}}$ and $\bar t=\sum_{i=1}^n T_{b_i}/n$, respectively. Both empirical quantities are based on $n=1000 000$ simulated FPTs, ensuring the closeness to the theoretical counterparts by the law of large numbers. This first scenario is meant to understand the goodness of our approximation through simulations.

Another relevant question is the performance of the MLEs and MEs of $\mu$ and $\sigma^2$, as described in Section \ref{Sec4b}. To compare different estimators, we use the relative mean error $R_\textrm{ME}$ to evaluate the bias and the relative mean square error $R_{\textrm{MSE}}$, which incorporates both the variance and the bias. They are defined as the average over the $1000$ repetitions of the quantities
\[
E_\textrm{rel}(\hat\mu)=\frac{\hat\mu-\mu}{\mu}, \qquad E_\textrm{rel sq}(\hat\mu)=\frac{(\hat\mu-\mu)^2}{\mu^2},
\]
and likewise for $\sigma^2$.
\begin{figure}
\includegraphics[width=1.0\textwidth]{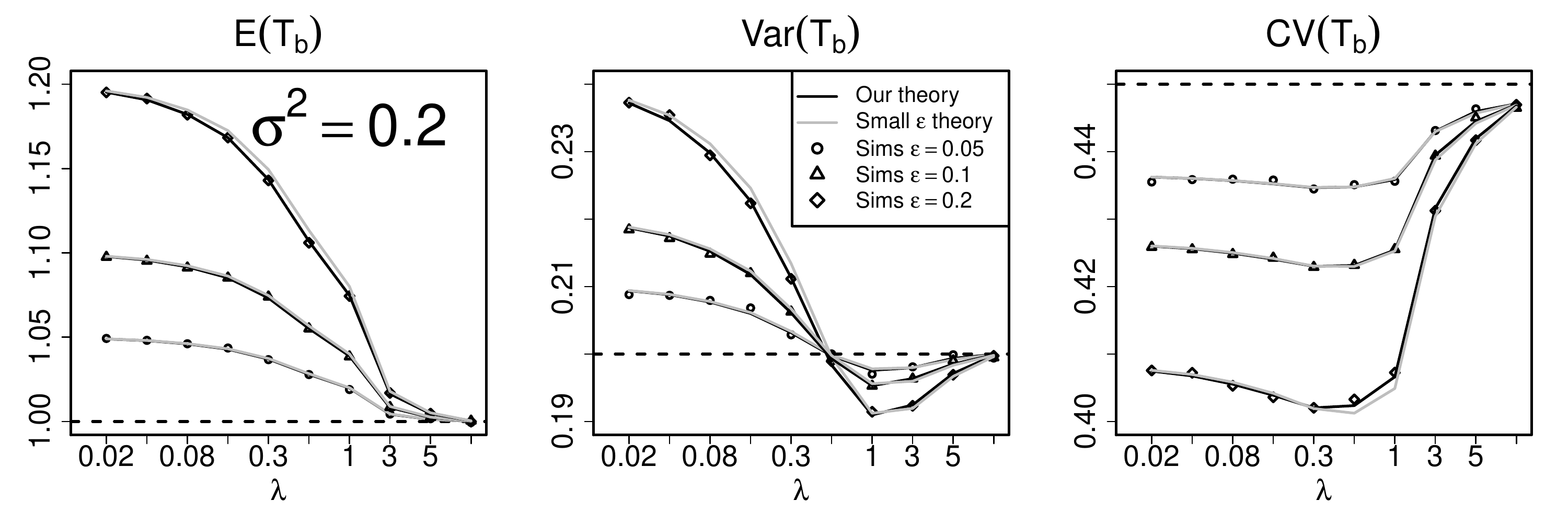}\\
\includegraphics[width=1.0\textwidth]{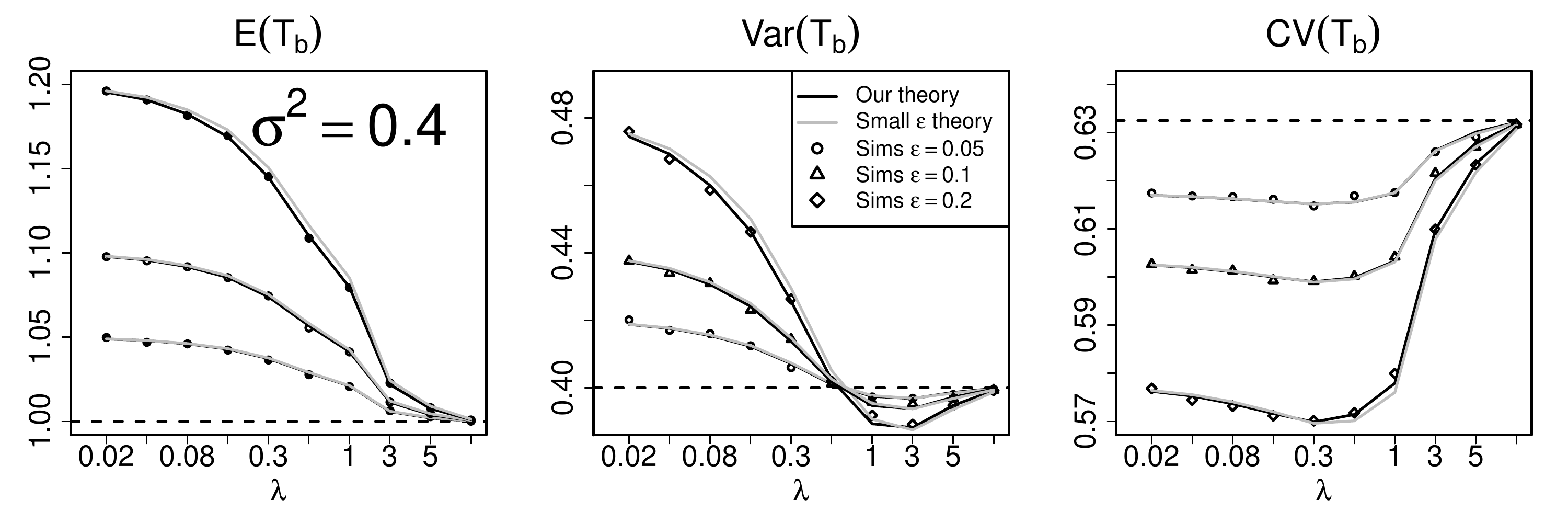}\\
\includegraphics[width=1.0\textwidth]{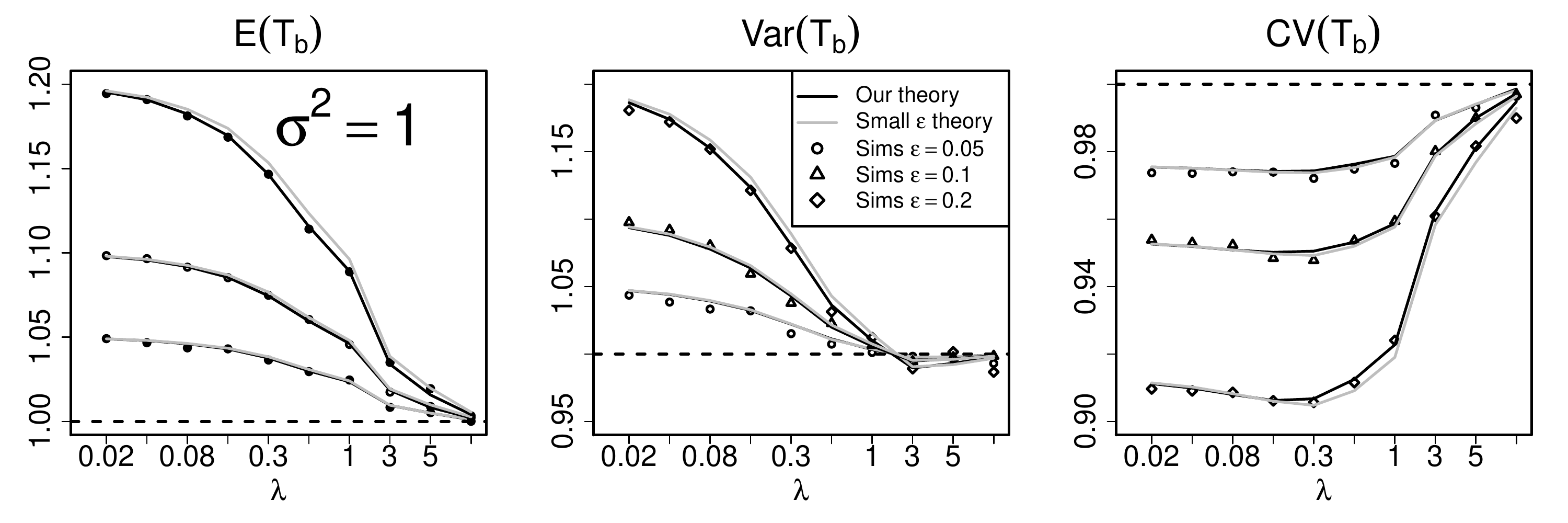}
\caption{Mean (left panels), variance (central panels) and CV (right panels) of the FPT $T_b$
 as a function of the decay rate $\lambda$ of the threshold for small values of the amplitude $\epsilon$ when $\mu=1$. Top panels: $\sigma^2=0.2$. Central panels: $\sigma^2=0.4$. Bottom panels: $\sigma^2=1$. Empirical quantities from simulations (symbols), theoretical quantities given by \eqref{EVT} for the piecewise linear threshold $b_\textrm{free}$ (solid lines), and theoretical quantities \eqref{ET} and \eqref{VT} when $\epsilon$ is small (solid gray lines). Also shown are the firing statistics of $T_b$ when $\epsilon=0$ (horizontal dashed lines).}
\label{Figall}
\end{figure}

\begin{figure}
\includegraphics[width=1.0\textwidth]{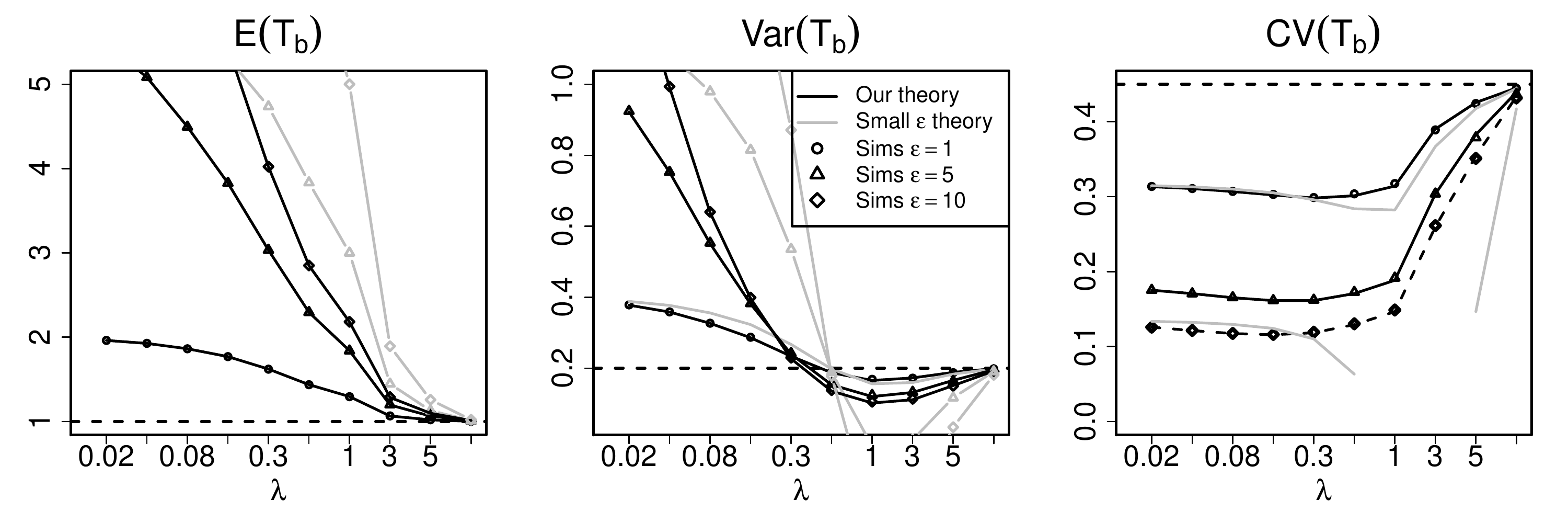}
\caption{Mean (left panel), variance (central panel) and CV (right panel) of the FPT $T_b$
 as a function of the decay rate $\lambda$ of the threshold for large values of the amplitude $\epsilon$ when $\mu=1$ and $\sigma^2=0.2$. Empirical quantities from simulations (symbols), theoretical quantities given by \eqref{EVT} for the piecewise linear threshold $b_\textrm{free}$ (solid lines), and theoretical quantities \eqref{ET} and \eqref{VT} when $\epsilon$ is small (solid gray lines). Also shown are the firing statistics of $T_b$ when $\epsilon=0$ (horizontal dashed lines).}
\label{Figall2}
\end{figure}

\subsection{Theoretical results for cdf and firing statistics of $T_b$} 
\begin{figure}
\includegraphics[width=\textwidth]{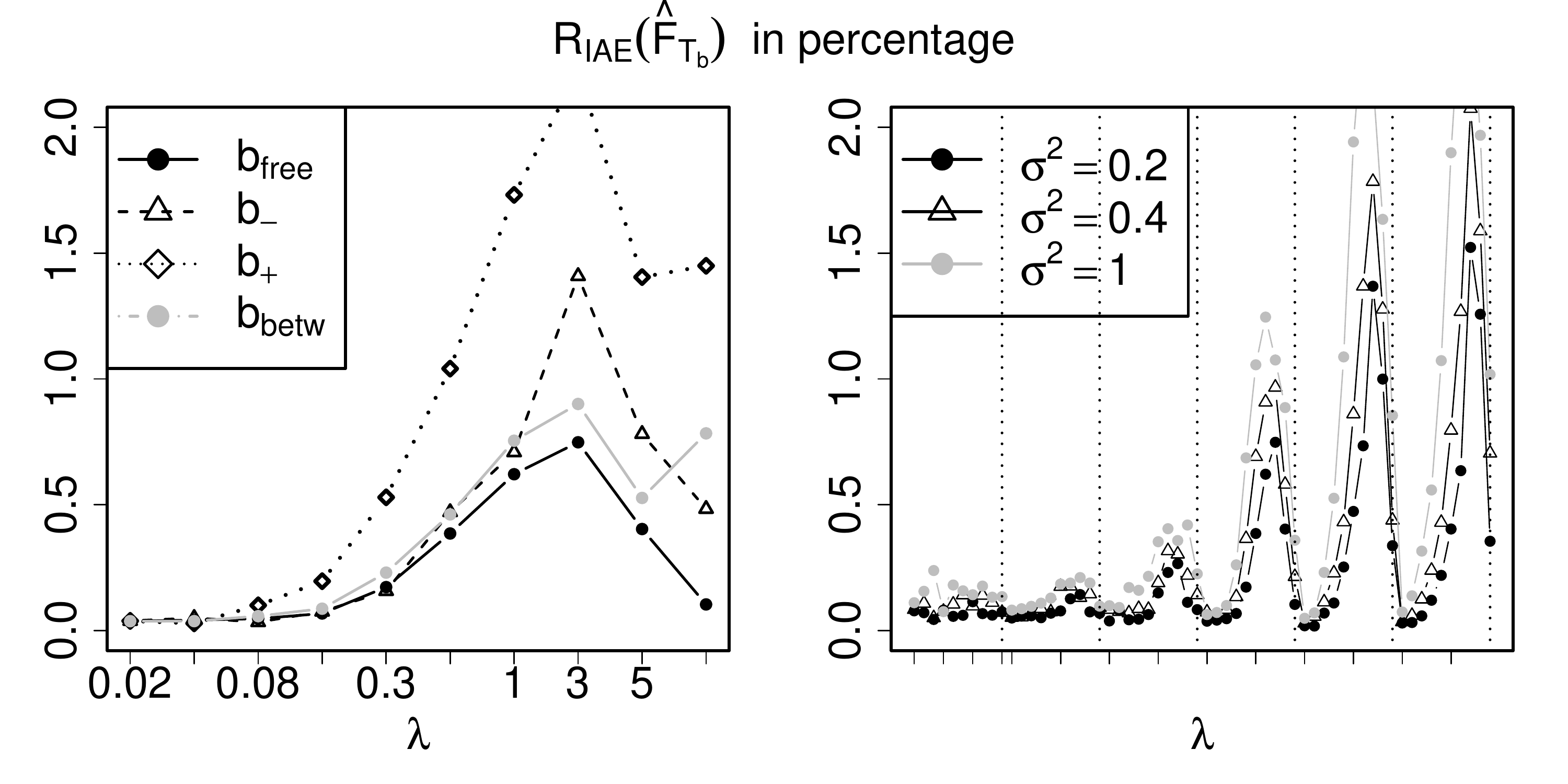}
\caption{$R_{\textrm{IAE}}(\hat F_{T_b})$ (in percentage) given by \eqref{riae} for different values of $\lambda$ and $\epsilon$ when $\mu=1$. Left panel: $R_{\textrm{IAE}}(\hat F_{T_b})$ from threshold $b_\textrm{free}$ (circles), $b_-$ (triangles), $b_+$ (rhombuses) and $b_\textrm{betw}$ (gray circles) when $\epsilon=1$ and $\sigma^2=0.2$.  Right panel: $R_\textrm{IAE}(\hat F_{T_{b_\textrm{free}}})$ for $\sigma^2=0.2$ (circles), $\sigma^2=0.4$ (triangles) and $\sigma^2=1$ (gray circles). The values of $\epsilon$ between consecutive vertical dotted lines are fixed and equal to $0.05, 0.1, 0.2, 1, 5,10$, while $\lambda$ varies between $0.02$ and $10$.}
\label{Figriae}
\end{figure}

In Fig. \ref{Figall} are reported theoretical and empirical means, variances and CVs of $T_b$ as a function of the rate $\lambda$, for small values of the amplitude $\epsilon$ and for $\sigma^2=0.2, 0.4$ and $1$.  The given theoretical quantities are obtained from \eqref{EVT} for the piecewise linear threshold $b_\textrm{free}$. Note how the mean of $T_b$ does not depend on $\sigma^2$, as it also happens for a linear threshold, while both variance and CV increase with growing $\sigma^2$. We refer to \cite{LindnerLongtin} for a detailed discussion on other qualitative features of the firing statistics, e.g. monotonic decrease on the mean with growing $\lambda$, existence of a minimum value for the variance, etc. What is relevant to emphasize is the excellent fit of the firing statistics provided by our method for any $\lambda$, and for both small (cf. Fig. \ref{Figall}) and large (cf. Fig. \ref{Figall2}) values of $\epsilon$. When $\epsilon$ is small, our theoretical firing statistics are at least as good as those in \cite{LindnerLongtin, Urdapilleta}, outperforming them when $\epsilon$ grows. The firing statistics of $T_{b_\textrm{betw}}$ are almost identical to those of $T_{b_\textrm{free}}$, while those of $T_{b_+}$ and $T_{b_-}$ are slightly different for increasing $\epsilon$. This can be seen in Fig. \ref{Figriae}, left panel, where the percentages of the $R_\textrm{IAE}(\hat F_T)$ for the four proposed estimators are given. As expected, the best approximation of $F_{T_b}$ is provided by $F_{T_{b_\textrm{free}}}$, since $b_\textrm{free}$ is the only threshold whose parameters are obtained from a non-constrained optimization problem. The performance of the estimators gets worse for large $\sigma^2$ and $\epsilon$. The highest error is observed for the value of $\lambda$ that minimizes the variance of $T_b$. However, all errors are smaller than $2\%$, confirming the good performance of the proposed estimators.

\subsection{Parameter estimation of $(\mu,\sigma^2)$}
\begin{figure}
\includegraphics[width=\textwidth]{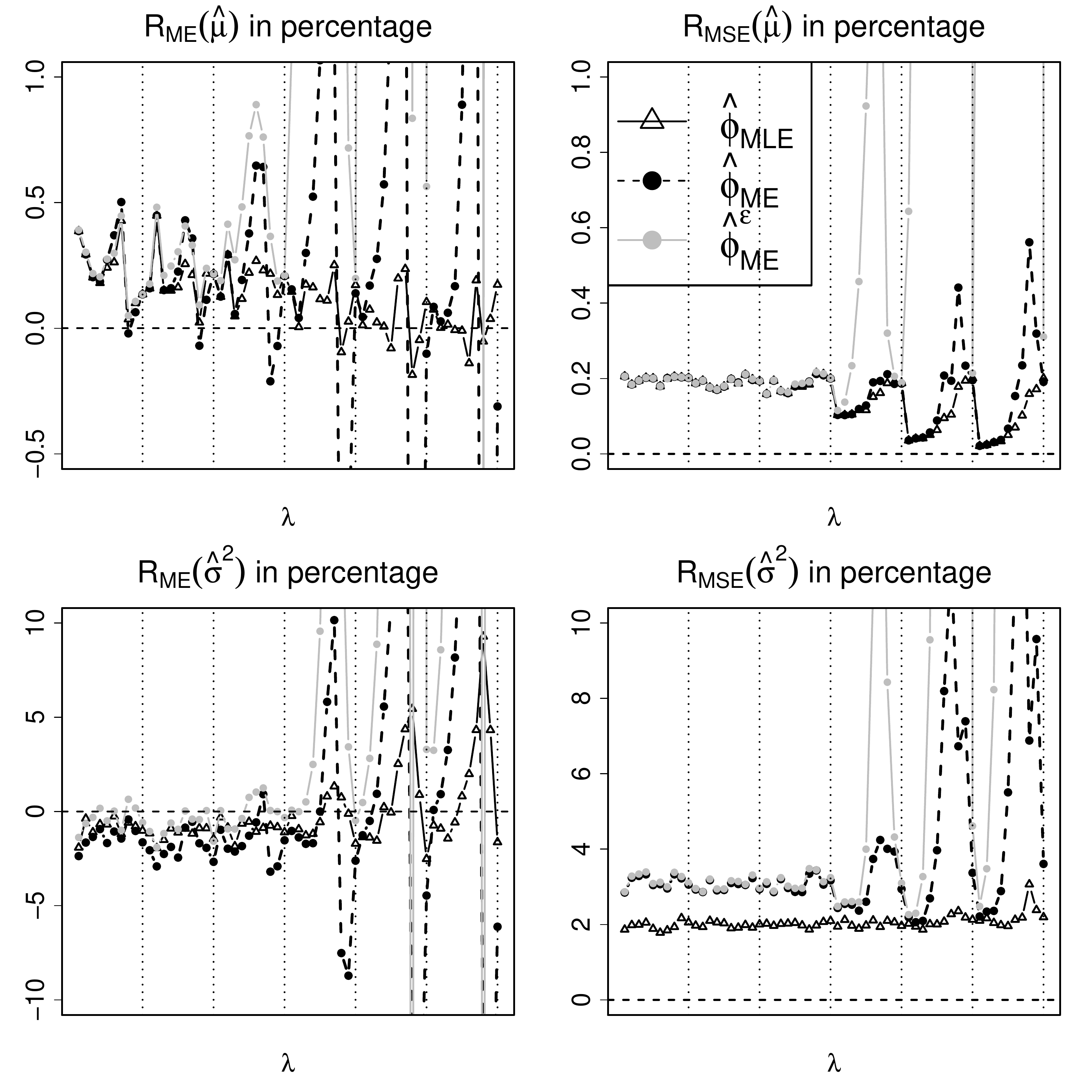}
\caption{Dependence of $R_\textrm{ME}(\hat\mu), R_\textrm{MSE}(\hat\mu), R_\textrm{ME}(\hat\sigma^2)$ and $R_\textrm{MSE}(\hat\sigma^2)$ (average over $1000$ simulations) on $\lambda$ and $\epsilon$ when $X$ is a Wiener process with $\mu=1$ and $\sigma^2=0.2$. Different estimators of $\phi=(\mu,\sigma^2)$ are considered: maximum likelihood estimator $\hat\phi_\textrm{MLE}$ (solid lines with triangles), moment estimator $\hat\phi_\textrm{ME}$ (dashed lines with circles) and moment estimator from \eqref{ET} and \eqref{VT} when $\epsilon$ is small, $\hat\phi_{ME}^\epsilon$ (gray solid lines with gray circles). The values of $\epsilon$ between consecutive vertical dotted lines are fixed and equal to $0.05, 0.1, 0.2, 1, 5,10$, while $\lambda$ varies between $0.02$ and $10$.
}
\label{Figphi1}
\end{figure}

\begin{figure}
\includegraphics[width=\textwidth]{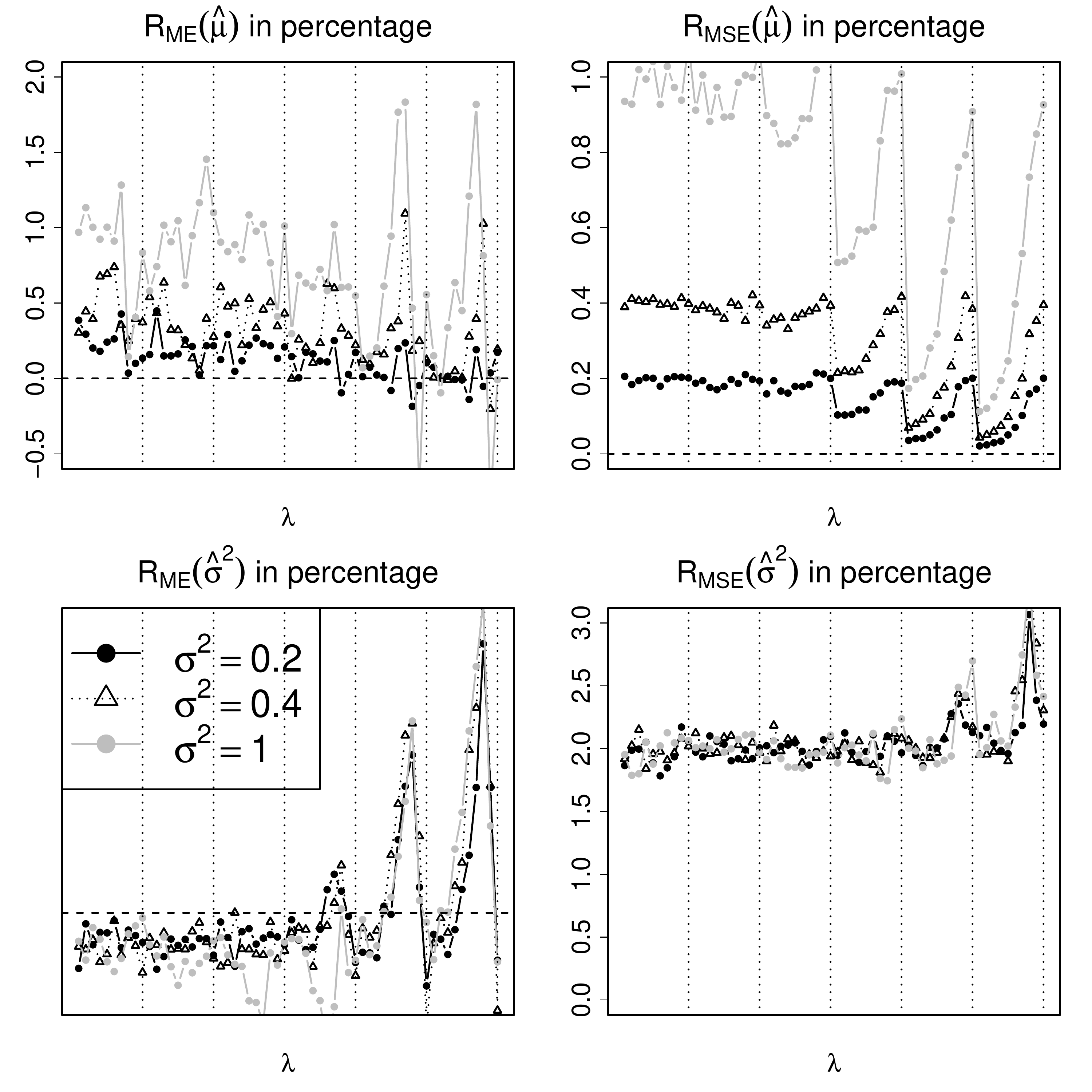}
\caption{Dependence of $R_\textrm{ME}(\hat\mu), R_\textrm{MSE}(\hat\mu), R_\textrm{ME}(\hat\sigma^2)$ and $R_\textrm{MSE}(\hat\sigma^2)$ (average over $1000$ simulations) on $\lambda, \epsilon$ and $\sigma^2$ when $X$ is a Wiener process with $\mu=1$ and $\sigma^2$ equal to $0.2$ (solid lines with circles), $0.4$ (dashed lines with triangles) and $1$ (gray solid lines with gray circles). Here only the maximum likelihood estimator $\hat\phi_\textrm{MLE}$ of $\phi=(\mu,\sigma^2)$ is considered. The values of $\epsilon$ between consecutive vertical dotted lines are fixed and equal to $0.05, 0.1, 0.2, 1, 5,10$, while $\lambda$ varies between $0.02$ and $10$.}
\label{Figphi2}
\end{figure}

We have seen that $T_{b_\textrm{free}}$ yields the best approximation of $T_b$ in terms of both cdf and firing statistics. For this reason, we limit our study to the estimators $\hat\phi$ based on $b_\textrm{free}$. In Fig. \ref{Figphi1} the $R_\textrm{ME}$ and the $R_\textrm{MSE}$ of $\hat\mu$ and $\hat\sigma^2$ are reported. As expected, the MLE provides the best estimate of $\phi$, while both MEs are acceptable only for small values of $\epsilon$. The performance of $\hat\mu$ is highly satisfactory, with $R_\textrm{ME}(\hat\mu)$ smaller than $0.5\%$, and $R_\textrm{MSE}(\hat\mu)<0.2\%$. Larger but still good $R_\textrm{ME}$ and $R_\textrm{MSE}$ are observed for $\hat\sigma^2$. The performance of $\hat\phi_\textrm{MLE}$ gets worse for growing $\sigma^2$, as shown in Fig. \ref{Figphi2}. However, except the $R_\textrm{ME}(\hat\sigma^2)$ for large values of $\epsilon$, all errors are between $0$ and $2-3\%$. Two last remarks should be done: first, the $R_\textrm{MSE}$ of $\hat\mu$ for small values of $\epsilon$ approaches the corresponding values of $\sigma^2$.
Second, $R_\textrm{MSE}(\hat\sigma^2)$ seems not to depend on $\lambda, \epsilon$ and $\sigma^2$, but to be equal to $2\%$. This error decreases when increasing the sample size. For example, the $R_\textrm{MSE}(\hat\sigma^2)\approx 1\%$ when $n=200$ (results not shown).

\section{Discussion} As a consequence of the recent increasing interest towards adapting-threshold models for the description of the neuronal spiking activity, a need of suitable mathematical tools to deal with hitting times of diffusion processes to time-varying thresholds arises. The mathematical literature on the FPT problem is rich and extensive. Unfortunately, analytical solutions are not available even for a problem as simple (compared to others) as the one considered here, i.e.  Wiener process to an exponentially decaying threshold. The closest result in this direction is represented by the work of Wang and P\"{o}tzelberger, who provide an explicit expression which should then be evaluated through Monte-Carlo simulations. The idea behind the works of Lindner and Longtin and of Urdapilleta, was to simplify some mathematical difficult equations arising from the study of the FPT by linearizing them in $\epsilon$, the amplitude of the decaying threshold. As a consequence, the quality of the approximation rapidly decreases when $\epsilon$ increases.

The method proposed here has no restriction on the parameter of the thresholds and it is based on the simple idea of replacing the boundary by a continuous two-piecewise linear threshold. This allows us to derive the analytical expression of the FPT density to the two-piecewise threshold, and to use it to approximate the desired distribution. To some extent, the presence of two linear thresholds can be considered as a second order approximation of the problem.

Numerical simulations show a good performance of the proposed method both when computing the main firing statistics, such as means, variances and CVs, and when calculating the FPT distribution. Different approximating thresholds have been proposed. We suggest choosing the one minimizing the distance with the curvilinear threshold and to restrict the interval where to perform the optimization as described in the paper. Among the estimators of the drift and diffusion coefficients of the Wiener process, we suggest applying MLE which always estimates the parameters  reasonably well.

The method proposed here may yield several interesting developments. First of all, it can be used to characterize the firing statistics of the Wiener process to the exponential decaying threshold, extending the previous considerations obtained for small values of $\epsilon$. Then, it may be extended to the case of a Wiener process with an adapting decaying threshold, as suggested in the introduction. Finally, our results may also be applied to all those processes which can be expressed as a piecewise monotone functional of a standard Brownian motion \cite{Wang2007}, as well as to Wiener processes with time-varying drift \cite{LindnerLongtin,Molini2011}.


\medskip
Received  June 02, 2015;  Accepted August 12, 2015.
\medskip

\end{document}